\documentclass[a4paper, 12pt,oneside,reqno]{amsart}
\usepackage[a4paper]{geometry}
\usepackage{amssymb,amsmath,amsthm}
\usepackage[T2A]{fontenc}
\usepackage{hyperref}

\usepackage{xcolor}

\def\Var{\mathrm{Var}}
\def\Nov{\mathrm{Nov}}
\def\Com{\mathrm{Com}}
\def\As{\mathrm{As}}
\def\GD{\mathrm{GD}}
\def\Lie{\mathrm{Lie}}
\def\Perm{\mathrm{Perm}}
\def\ComTriAs{\mathrm{ComTriAs}}
\def\Zinb{\mathrm{Zinb}}
\def\Leib{\mathrm{Leib}}
\def\Pois{\mathrm{Pois}}

\def\pre{\mathrm{pre}}
\def\post{\mathrm{post}}
\def\di{\mathrm{di}}
\def\tri{\mathrm{tri}}
\def\NP{\mathrm{NP}}
\def\Res {\mathop {\fam 0 Res}}

\def\End {\mathop {\fam 0 End}\nolimits}
\def\Cend {\mathop {\fam 0 Cend}\nolimits}
\def\Hom {\mathop {\fam 0 Hom}\nolimits}

\def\vecotimes{\mathbin{\vec\otimes}}
\def\oo#1{\mathbin{{}_{(#1)}}}
\def\ooc#1#2{\mathbin{{\circ}_{#1(#2)}}}

\newtheorem{definition}{Definition}
\newtheorem{lemma}{Lemma}
\newtheorem{proposition}{Proposition}
\newtheorem{theorem}{Theorem}
\newtheorem{corollary}{Corollary}
\newtheorem{remark}{Remark}
\newtheorem{example}{Example}

\title{On the Dong Property for a binary quadratic operad}
\author{P.~S.~Kolesnikov, B.~K.~Sartayev}
\address{Sobolev Institute of Mathematics, Novosibirsk (Russia)}
\address{Narxoz University, Almaty (Kazakhstan)}

\begin{document}

\begin{abstract}
The classical Dong Lemma for distributions over a Lie algebra lies in the foundation of vertex algebras theory.
In this paper, we find necessary and sufficient condition for a variety of nonassociative algebras with binary operations to satisfy the analogue of the Dong Lemma.
In particular, it turns out that Novikov and Novikov--Poisson algebras satisfy the Dong Lemma. The criterion is stated in the language of operads, so we determine for which binary quadratic operads the Dong Lemma holds true.
As an application, we show the black Manin product of Dong operads is also a Dong operad.
\end{abstract}

\maketitle

\section{Introduction}
One of the cornerstones of the theory of vertex algebras 
is the statement known as the Dong Lemma for quantum fields
(see, e.g., \cite{FbZ2003}) which goes back to the paper 
\cite{HLie199x}. The Dong Lemma concerns the locality property 
of formal distributions with coefficients in a Lie algebra, but 
the roots of this statement are in axiomatic quantum field theory in the sense of \cite{Wightman5x}: the locality of formal distributions 
represents the basic principle of relativity transformed 
into algebraic conditions in the special case (2-dimensional conformal field theory) \cite{BPZ198x}. 

The first algebraic definition of a vertex algebra is due to \cite{Bor198x}, and a review of modern approaches to the subject may
be found, for example, in \cite{BK-et-al2018}. The Dong Lemma is 
especially important for constructing examples of 
vertex algebras. Namely, 
given a family of pairwise local quantum fields, 
we are sure that the locality property holds
for the entire space of quantum fields generated (in the context of vertex algebras) by the initial ones.

In particular, the Dong Lemma plays a crucial role in the theory 
of conformal algebras. The latter were introduced in \cite{KacVABeginn} as an algebraic tool to describe the singular 
part of a vertex algebra operations. Namely, suppose
$(V, T, |0\rangle , Y)$ is a vertex algebra with a space of states
$V$, a translation operator $T: V\to V$, a vacuum vector $|0\rangle \in V$, and vertex operator 
$Y: V\to \mathrm{gl}\,(V)[[z,z^{-1}]]$,
$a\mapsto Y(a,z)$.
One of the most important conditions on a vertex operator
is called {\em locality}:
for every $a,b\in V$ there exists $N\in \mathbb Z_+$
such that 
\[
(w-z)^N [Y(a,w),Y(b,z)] = 0
\]
is the space $\mathrm{gl}\,(V)[[z,z^{-1},w,w^{-1}]]$.

Then $V$ equipped with the same operator $T$ and with a series of 
$n$-{\em products} $[\cdot \oo{n}\cdot ]$, $n\in \mathbb Z_+$,
such that
\begin{equation}\label{eq:n-prodY}
Y([a\oo{n} b], z) = \Res\limits_{w=0} 
 (w-z)^n[Y(a,w),Y(b,z)],
\quad a,b\in V,
\end{equation}
is a {\em Lie conformal algebra} \cite{KacVABeginn}. 

Thus, the notion of a Lie conformal algebra describes a part of
a vertex algebra structure. However, for every Lie conformal algebra 
there exists a uniquely defined universal enveloping vertex algebra, and most of classical examples of vertex algebras (see, e.g., \cite{FbZ2003}) are based on this construction.

An axiomatic description of (Lie and associative) conformal algebras is given in \cite{KacVABeginn}.
The latter notion naturally emerges in the study of representations of Lie conformal algebras. In recent years, a series of papers has been devoted to conformal algebras corresponding to various classes of nonassociative algebras (e.g., Jordan \cite{KacRetakh-JordSuper}, left-symmetric \cite{Hong2015}, dendriform \cite{Yuan22}). A general method to get a unified approach to all these definitions was described in \cite{BDK2001}: 
given an operad $\mathcal P$, the class of 
$\mathcal P$-conformal algebras consists of all morphisms from 
$\mathcal P$ to the symmetric multi-category (pseudo-tensor category) $\mathcal M^*(H)$, where $H=\mathbb C[T]$ is the polynomial Hopf algebra in one primitive variable. 

For every $\mathcal P$-conformal algebra $V$ there exists an 
``ordinary'' $\mathcal P$-algebra $A=\mathcal A(V)$ such that 
$V$ can be embedded into the space of formal distributions over $A$, i.e., there exists an injection
\[
V\ni a \mapsto Y(a,z) \in A[[z,z^{-1}]],
\]
so that the operator $T$ corresponds 
to the formal derivation $d/dz$ and the conformal $n$-products 
are given by expressions similar to \eqref{eq:n-prodY}.

In order to construct $\mathcal P$-conformal algebras 
by means of formal distributions over ordinary 
$\mathcal P$-algebras, 
it is essential to know whether the locality of formal 
distributions is preserved under $n$-products.

Let $\mathfrak V$ be a class of algebras with a family 
of bilinear products 
$\mu_i$, $i\in I$ (later we will always assume $I$ is a finite set).
Suppose that for every $\mu_i$ the opposite operation $(12)\mu_i$ also belongs to this family.
(Practically, we consider a space of binary operations 
as a module over $S_2$, $(\mu_i)_{i\in I}$ is a basis of this space.) The operations $\mu_i$ are expanded on 
formal distributions $a(w),b(z)\in A[[z,z^{-1},w,w^{-1}]]$
in the ordinary way (formal variables commute 
with each other and with coefficients).

Given an algebra $A\in \mathfrak V$ and 
two formal distributions $a(z),b(z)\in A[[z,z^{-1}]]$, 
we say $(a(z),b(z))$ is a {\em local} pair 
if for every $i\in I$ there exists $N\in \mathbb Z_+$
such that 
\[
(w-z)^N \mu_i(a(w),b(z)) = 0.
\]

For a local pair $(a(z),b(z))$ of formal distributions 
over $A\in \mathfrak V$ define 
\begin{equation}\label{eq:n-prod}
(a\ooc{i}{n} b)(z) = \Res\limits_{w=0} 
 (w-z)^n \mu_i(a(w),b(z)),
\quad i\in I,\ n\in \mathbb Z_+.
\end{equation}
The formal series 
$(a\ooc{i}{n} b)(z)$ are the coefficients in the 
expansion of the product $\mu_i(a(w),b(z))$ 
into a finite sum
\[
\mu_i(a(w),b(z)) = \sum\limits_{n=0}^{N-1}
(a\ooc{i}{n} b)(z)\cdot \dfrac{1}{n!} \partial^n_z\delta(w-z), 
\]
where $\delta(w-z)$ is the formal delta-function 
\cite{KacVABeginn}.

We say $\mathfrak V$ {\em satisfies the Dong Property}
if for every $A\in \mathfrak V$ and
for every pairwise local formal distributions 
$a(z),b(z),c(z)\in A[[z,z^{-1}]]$
the pair $((a\ooc{i}{n} b)(z), c(z))$ is local
for all $i\in I$, $n\in \mathbb Z_+$.

If $\mathfrak V$ is a variety of algebras, i.e., 
a class defined by a family of identities, 
then the Dong Property obviously concerns 
multilinear identities of degree~3. 
Hence, it is natural to formulate the Dong Property 
for a binary quadratic operad. 
Namely, we say that an operad $\mathcal P$ satisfies 
the Dong Property (or $\mathcal P$ is a {\em Dong operad}) 
if so is the class of all 
$\mathcal P$-algebras.
Apparently, a binary operad $\mathcal P$ satisfies the 
Dong Property if and only if so is its quadratic part 
(see, e.g., \cite[S.\,7.8]{LodayVallette_book}).

In this paper, we establish the 
necessary and sufficient 
condition for a binary quadratic operad 
$\mathcal P$ to satisfy the Dong Property. 
The condition is easy to check provided that we know 
the Koszul dual operad $\mathcal P^!$. 

As an application, we observe how operads with the
Dong Property behave under some standard operadic constructions.
For example, 
we show that for neither binary quadratic operad $\mathcal P$ its dendriform splitting $\pre\mathcal P$ or $\post\mathcal P$ satisfies the Dong Property.
If $\mathcal P$ is a Dong operad then so are 
replicated (\cite{BGG_HK_Proceed}) operads 
$\di\mathcal P$ and $\tri\mathcal P$. 
The latter are obtained by means of Manin white product 
of $\mathcal P$ with operads $\Perm$ and $\ComTriAs$, 
respectively. However, the white Manin product of two 
Dong operads is not necessarily a Dong operad 
(it is well known, for example, that the white product 
of $\As$ and $\Lie $ corresponds to the class of all nonassociative algebras).
Finally, we consider the black Manin product of two Dong operads and show that it is always a Dong operad. 
To that end, we establish an easy algorithmic way to compute the Manin black product. 

Throughout the paper, 
all spaces are over a field $\Bbbk $,
$S_n$ is the symmetric group. 
Given two $S_n$-modules $U$ and $V$ we denote by 
$U\otimes V$ the ordinary tensor product with 
$\sigma (u\otimes v) = \sigma u\otimes \sigma v$, 
$\sigma \in S_n$, $u\in U$, $v\in V$.
We will use the notation $U\vecotimes V$ 
for the right-justified $S_n$-module:
$\sigma (u\vecotimes v) = u\vecotimes \sigma v$.
The dual space $V^*$ equipped with the transpose dual action
of $S_n$ is denoted by $V^\vee $. Obviously, 
$V^\vee \cong \Bbbk_-\otimes V^*$, 
where $\Bbbk_-$ is the 1-dimensional sign $S_n$-module, 
$(U\vecotimes V)^\vee \cong U^\vee\vecotimes V^\vee$, 
$(U\otimes V)^\vee \cong \Bbbk_-\otimes U^\vee\otimes V^\vee$. For the induced modules ($m\ge n$) we have 
\[
\mathrm{Ind}_{S_n}^{S_m} (V^\vee ) \cong \big(
\mathrm{Ind}_{S_n}^{S_m} (V) 
\big)^\vee .
\]

\section{Binary quadratic operads and Koszul duality}

An operad $\mathcal O$ in our setting is a family of 
linear spaces $\mathcal O(n)$, $n=1,2,3,\ldots$, such that 
every $\mathcal O(n)$ is a (left) $S_n$-module, and the 
following composition rules are given:
\[
\gamma^n_{m_1,\ldots, m_k} : 
\mathcal \mathcal O(k)\otimes O(m_1)\otimes \dots \otimes \mathcal O(m_k)
\to \mathcal O(n),\quad n=m_1+\dots + m_k.
\]
The composition rules are linear maps such that 
the natural associativity and equivariance (relative to the action of $S_n$) properties hold
(see, e.g., \cite{LodayVallette_book, BrD_Companion} 
for details).
The space $\mathcal O(1)$ contains an identity element $1_{\mathcal O}$ with respect to the composition.

\begin{example}
Given a linear space $A$, define 
an operad $\End_A$ in such a way that 
$\End_A(n) = \Hom (A^{\otimes n}, A)$. 
This is an operad relative to the ordinary 
composition of multilinear maps 
(the groups $S_n$ act by permutation 
of arguments).
\end{example}

\begin{example}
Let $H$ be a cocommutative bialgebra 
with a coproduct $\Delta : H\to H\otimes H$. 
Then for every (left, unital) $H$-module $M$
one may define 
an operad $\Cend_M$ in such a way that 
\[
\Cend_M(n) = \Hom_{H^{\otimes n}}(M^{\otimes n}, 
H^{\otimes n}\otimes_H M), 
\]
see \cite{BDK2001} or \cite{Kol2006CA} for details.
\end{example}

\begin{example}\label{exmp:VarOperad}
Let $\Var $ be a multilinear variety of algebras, 
i.e., a class of linear spaces equipped 
with a series of multilinear operations 
(not necessarily binary) satisfying a family 
of mutlilinear identities in the variables 
$X=\{x_1,x_2,\dots \}$.
Then the multilinear part of the free algebra 
$\Var \langle X\rangle $ generated by $X$ is an operad (also denoted $\Var$) so that 
$\Var (n)$, $n\ge 1$, consists of all multilinear 
elements of degree $n$ depending on $x_1,\ldots, x_n$. Symmetric groups $S_n$ permute the variables, 
and the composition rules are constructed by substitutions with consecutive renumeration. 
The identity element in $\Var(1)$ is always $1_{\Var }= x_1$.
\end{example}

For example, let $\Var =\As$ be the variety 
of associative algebras. Then 
$\As(n)$ is spanned by $x_{\sigma(1)}\dots x_{\sigma(n)}$, $\sigma \in S_n$, and
\[
\gamma_{2,1,3}^6: 
x_3x_2x_1\otimes 
x_2x_1\otimes x_1\otimes x_1x_3x_2
\mapsto 
(x_4x_6x_5)x_3(x_2x_1) = x_4x_6x_5x_3x_2x_1.
\]

A morphism of operads $\mathcal P\to \mathcal O$
is a family of $S_n$-linear maps 
$\mathcal P(n)\to \mathcal O(n)$
that preserve compositions and identity.
For example, if $\mathcal O=\As$ as above 
and $\mathcal P=\Lie$ is the operad obtained from the variety of Lie algebras 
then there exists a morphism of operads 
$\Lie \to \As$
transforming $x_1x_2\in \Lie(2)$ into $x_1x_2-x_2x_1\in \As(2)$. 

A binary quadratic operad is a combinatorial object which is 
convenient to use for the description of a class of algebras 
with (one or more) bilinear operations satisfying multilinear 
identities of degree~3. Let us state the necessary details, mainly following \cite{GinKapr1994} in order to fix the notations.

Let $V$ be a finite-dimensional space over a field $\Bbbk $ 
equipped with a linear action of the symmetric group $S_2$.
Define 
\[
\mathcal F_V(3) = \Bbbk S_3\otimes _{\Bbbk S_2} (V\vecotimes V).
\]
Given an $S_3$-submodule $R$ in $\mathcal F_V(3)$, one can define 
a binary quadratic operad $\mathcal P = \mathcal P(V,R)$ so that 
$\mathcal P(0)=\Bbbk $, $\mathcal P(2)=V$, 
$\mathcal P(3) = \mathcal F_V(3)/R$.

Given a binary quadratic operad $\mathcal P=\mathcal P(V,R)$ and a linear space 
$A$ over the field $\Bbbk $, a morphism of operads 
\[
\mathcal P \to \End_A
\]
defines an algebra structure on the space $A$ in such a way 
that the image of $V=\mathcal P(2)$ is the 
subspace of $\End_A(2)$ spanned by the binary operations in $A$ (together with their opposites). The class of $\mathcal P$-algebras obtained 
in this way is a multilinear 
variety with operations $V$ (for a more classical setting, one may fix a linear basis of $V$) and identities given by $R$.

It is worth mentioning that the same idea 
underlies the categorial definition of an $H$-conformal 
$\mathcal P$-algebra: a morphism of operads 
\[
\mathcal P\to \Cend_M
\]
defines an  $H$-conformal algebra (or pseudo-algebra \cite{BDK2001}) structure on a module $M$
over a cocommutative bialgebra~$H$.
The case when $H=\Bbbk [\partial ]$ with 
$\Delta(\partial)=\partial\otimes 1+1\otimes \partial $, $\varepsilon(\partial )=0$
corresponds to what is known as conformal algebras
\cite{KacVABeginn}.

The interpretation of the space $\mathcal P(3)$
in the language of Example~\ref{exmp:VarOperad}
is the following. 
If $e_i$ and $e_j$ are from $V$, and 
$e_i=x_1\circ_i x_2$, $e_j = x_1\circ_j x_2$
then 
$e_i\vecotimes e_j =(x_1\circ_j x_2)\circ_i x_3$.

For example, if $\dim V=2$ and a basis of $V$ is presented by $\mu=x_1x_2$, $(12)\mu=x_2x_1$ 
then the space $V\vecotimes V$ is 4-dimensional, 
it is spanned by
\[
\mu\vecotimes \mu = (x_1x_2)x_3,
\quad 
\mu\vecotimes (12)\mu = (x_2x_1)x_3,
\quad
(12)\mu\vecotimes \mu = x_3(x_1x_2),
\quad
(12)\mu\vecotimes (12)\mu = x_3(x_2x_1).
\]

\begin{example}\label{exmp:Pois}
The operad $\Pois $ governing the 
variety of Poisson algebras
is defined by the following pair $(V,R)$. 
The space $V$ is a 2-dimensional $S_2$-module 
with a basis 
$(e_1,e_2)$, where $(12)e_1=e_1$, $(12)e_2=-e_2$.
The $S_3$-submodule $R$ of $\mathcal F_V(3)$ 
generated by
\begin{equation}\label{eq:OperadicAssoc-e_1}
e_1\vecotimes e_1 - (13)(e_1\vecotimes e_1)
\end{equation}
(the associativity of $e_1=x_1x_2$), 
\begin{equation}\label{eq:OperadicJacobi-e_2}
e_2\vecotimes e_2 - (13)(e_2\vecotimes e_2)
- (23)(e_2\vecotimes e_2)
\end{equation}
(the Jacobi identity for $e_2=\{x_1,x_2\}$),
and 
\[
e_2\vecotimes e_1 -(23)(e_1\vecotimes e_2)
 +(13)(e_1\vecotimes e_2)
\]
(the Leibniz rule for $\{x_1x_2,x_3\}$).
\end{example}

Alternatively, one may consider another basis 
of $\Pois(2)$ that consists of $\mu=e_1+e_2$
and $(12)\mu = e_1-e_2$. Then it is enough to set one identity to define the same $R$
(see \cite{MarklRemm}).

\begin{example}\label{exmp:Nov}
The operad $\Nov $ of Novikov algebras
is defined by a 2-dimensional space $V$
with a basis $g_1, g_2=(12)g_1$ and by 
the $S_3$-submodule $R$ generated by
\begin{equation}\label{eq:OperadicLSym-g_1g_2}
g_1\vecotimes g_1 
-g_1\vecotimes g_2 
-(13)(g_2\vecotimes g_2)
+(23)(g_2\vecotimes g_1 )
\end{equation}
(the left symmetry of $g_1=x_1\circ x_2$),
\begin{equation}\label{eq:OperadicRCom-g_1g_2}
g_1\vecotimes g_1 - (23)(g_1\vecotimes g_1)
\end{equation}
(the right commutativity for $x_1\circ x_2$).
\end{example}

The class of Novikov algebras was introduced in the 1980s as a language to describe the series of conditions on the coefficients of a rank 3 tensor emerging in the formal calculus of variations \cite{GelDorfm} or in the construction
of the Hamiltonian formalism for PDEs of hydrodynamic type \cite{BalNov}. 

In \cite{Xu1993}, the notion of a Novikov--Poisson 
algebras (NP-algebras) was introduced as a tool to study Novikov algebras. When translated into the language of binary quadratic operads, this definition
turns into the following

\begin{example}\label{exmp:NovPois}
The operad $\NP$ is defined by a 3-dimensional 
space $V$ with a basis $e_1,g_1,g_2$, where 
$(12)e_1=e_1$, $(12)g_1=g_2$, 
relative to the $S_3$-submodule $R$ generated by
\eqref{eq:OperadicAssoc-e_1},
\eqref{eq:OperadicLSym-g_1g_2},
\eqref{eq:OperadicRCom-g_1g_2},
and
\begin{equation}\label{eq:OperadicNP-e_1g_1g_2}
g_1\vecotimes e_1 - (13)(e_1\vecotimes g_2),
\quad 
e_1\vecotimes g_1 -(13)(g_2\vecotimes e_1)
-e_1\vecotimes g_2 +(23)(g_2\vecotimes e_1).
\end{equation}
Straightforward computation shows the dimension of $R$ as of a linear space is equal to~16.
\end{example}

The class of Novikov--Poisson algebra is somewhat analogous to the class of Poisson algebras with a Novikov product 
instead of the Lie bracket. 
The classes of Lie and Novikov algebras have many similar features. Note that the original name for the class of Novikov algebras was 
``local translation invariant Lie algebras'' \cite{Zelmanov87}, the name ``Novikov algebras'' was proposed in \cite{Osborn90}. In this paper, we will  mention one more common feature of Poisson and NP-algebras: they both satisfy the Dong Property.

Let us also mention an example of a binary quadratic operad corresponding to a variety of algebras 
introduced in \cite{GelDorfm}, see also \cite{Xu1999}.

\begin{example}\label{exmp:GD}
The operad $\GD$ of Gelfand--Dorfman algebras
is generated by 3-dimensional space $V$ with a basis
$g_1,g_2,e_2$, where 
$(12)g_1=g_2$, $(12)e_2=-e_2$, relative to the 
$S_3$-submodule $R$ generated by
\eqref{eq:OperadicJacobi-e_2},
\eqref{eq:OperadicLSym-g_1g_2},
\eqref{eq:OperadicRCom-g_1g_2}, 
and
\begin{equation}\label{eq:OperadicGD}
e_2\vecotimes g_1 - (23)(e_2\vecotimes g_1 )
+ g_1\vecotimes e_2 -(23)(g_1\vecotimes e_2)
+(13)(g_2\vecotimes e_2).
\end{equation}
\end{example}

Given a binary quadratic operad 
$\mathcal P=\mathcal P(V,R)$, 
its Koszul dual operad 
$\mathcal P^!$ is determined by 
the space of generators $V^\vee $ 
and relations 
$R^\perp \subseteq (\mathcal F_V(3))^\vee
\cong \mathcal F_{V^\vee}(3)$, 
see \cite{GinKapr1994}.
There is a relatively simple way to compute 
the defining identities of the variety of 
$\mathcal P^!$-algebras
without solving a large system of linear equations
describing~$R^\perp $.

Suppose $(e_i)_{i\in I}$ is a linear basis 
of $V$, the corresponding operations on 
a $\mathcal P$-algebra are denoted 
$(\cdot \circ_i\cdot)$, $i\in I$. 
Recall that $|I|$ is assumed to be finite.

Consider the dual basis 
$(e_i^\vee )_{i\in I}$ of $V^\vee $.
The corresponding operations 
on a $\mathcal P^!$-algebra 
are denoted $(\cdot \circ^i \cdot)$, $i\in I$.
Then the elements of 
$\mathcal F_{V^\vee}(3)$ are presented by 
monomials of the form 
$(x_{\sigma(1)}\circ^i x_{\sigma(2)})\circ^j x_{\sigma(3)}$, for $i,j\in I$, $\sigma \in S_3$.
The defining relations of the operad $\mathcal P^!$ can be found by the following statement.

\begin{proposition}[c.f. {\cite[2.2.9]{GinKapr1994}}]\label{prop:DualOp-Lie}
Let $\mathcal P=\mathcal P(V,R)$ be a binary quadratic operad,
and let $B$ be an algebra with
the $S_2$-space  $V^\vee$ of bilinear operations. 
Then $B$ is a $\mathcal P^!$-algebra 
if and only if 
for every $\mathcal P$-algebra $F$
the space 
$B\otimes F$ equipped with the bracket 
\begin{equation}\label{eq:DualProdLie}
 [(a\otimes x), (b\otimes y)]   
 =\sum\limits_{i\in I} (a\circ^i b)\otimes (x\circ_i y),
 \quad a,b\in B,\ x,y\in F,
\end{equation}
is a Lie algebra (i.e., $[\cdot,\cdot]$ satisfies the Jacobi identity).
\end{proposition}

\begin{example}\label{exmp:NP-dual}
Let $\mathcal P=\NP$ be the operad of 
Novikov--Poisson algebras from Example~\ref{exmp:NovPois}. 
Assume
$e_1 = x_1x_2$,  $g_2=x_1\circ x_2$,
$e_3 = x_2\circ x_1$.

Consider the dual basis 
$e_1^\vee, g_1^\vee, g_2^\vee $, and denote, 
respectively, 
$e_1^\vee = [y_1,y_2]$, 
$g_1^\vee = y_1\circ y_2$,
then 
$e_3^\vee = - y_2\circ y_1$.
By Proposition~\ref{prop:DualOp-Lie}, 
the relations on $e_1^\vee, g_1^\vee, g_2^\vee$
are exactly those that make 
the skew-symmetric bracket 
\[
[y_1\otimes x_1, y_2\otimes x_2]
= e_1^\vee \otimes e_1 + g_1^\vee\otimes g_1 + g_2^\vee \otimes g_2 
= [y_1,y_2]\otimes (x_1x_2) 
+ (y_1\circ y_2)\otimes (x_1\circ x_2)
- (y_2\circ y_1)\otimes (x_2\circ x_1)
\]
to satisfy the Jacobi identity. 
\end{example}

Calculate the first term of the Jacobi identity: \begin{multline*}
[[y_1\otimes x_1,y_2\otimes x_2],y_3\otimes x_3]
=[[y_1,y_2],y_3]\otimes (x_1x_2)x_3
+[y_1,y_2]\circ y_3 \otimes (x_1x_2)\circ x_3
-y_3\circ[y_1,y_2]\otimes x_3\circ(x_1x_2)
\\
+[y_1\circ y_2,y_3]\otimes (x_1\circ x_2)x_3
+(y_1\circ y_2)\circ y_3 \otimes (x_1\circ x_2)\circ x_3
-y_3\circ(y_1\circ y_2)\otimes x_3\circ(x_1\circ x_2)
\\
-[y_2\circ y_1,y_3]\otimes (x_2\circ x_1)x_3
+(y_2\circ y_1)\circ y_3\otimes (x_2\circ x_1)\circ x_3
-y_3\circ(y_2\circ y_1)\otimes x_3\circ(x_2\circ x_1).
\end{multline*}
Apply the cyclic permutations $(123)$, $(132)$ to get 
the remaining terms of the Jacobi identity, then add them all together and collect similar terms at normal forms (basic elements of $\NP(3)$) 
in the second tensor factor.
These normal forms are: 
$x_1x_2x_3$, 
$x_1x_3\circ x_2$,
$x_2\circ x_1x_3$,
$(x_1\circ x_3)x_2$, 
$(x_3\circ x_1)x_2$, 
$x_1\circ (x_2\circ x_3)$,
$x_1\circ (x_3\circ x_2)$,
$x_2\circ (x_1\circ x_3)$,
$x_2\circ (x_3\circ x_1)$,
$(x_1\circ x_3)\circ x_2$,
$(x_3\circ x_1)\circ x_2$.
As a result, we obtain the following identities
for $e_1^\vee = [y_1,y_2]$ (skew symmetric)
and $g_1^\vee = y_1\circ y_2$:
\begin{gather}
{}
[[y_1,y_2],y_3]+[[y_2,y_3],y_1]+[[y_3,y_1],y_2]=0,
  \nonumber 
\\
(y_1\circ y_2)\circ y_3-y_1\circ(y_2\circ y_3)
=(y_1\circ y_3)\circ y_2 - y_1\circ(y_3\circ y_2),
 \quad 
y_1\circ (y_2\circ y_3)=y_2\circ (y_1\circ y_3), 
\nonumber 
\\
y_1\circ[y_2,y_3]+y_2\circ[y_3,y_1]+y_3\circ[y_1,y_2]=0, \label{eq:NP-Dual-1} 
\\
{}
[y_1\circ y_2, y_3]+[y_1,y_3\circ y_2]-[y_1,y_3]\circ y_2-y_2\circ[y_1,y_3]=0.
  \label{eq:NP-Dual-2}
\end{gather}
One may easily see that \eqref{eq:NP-Dual-1},
\eqref{eq:NP-Dual-2} imply an 
identity which is opposite to \eqref{eq:OperadicGD}. 
Hence, quite expectedly, 
if $V$ is an $\NP^!$-algebra then 
the opposite algebra $V^{op}$ is a $\GD$-algebra. 

\begin{remark}
Note that every transpose Poisson algebra 
in the sense of \cite{BaiKo-TPalgebras}
is also an $\NP^!$-algebra.
\end{remark}

\section{The Dong Property: sufficient condition}

Let $\mathcal P=\mathcal P(V,R)$ be a binary 
quadratic operad (we assume $\dim V<\infty$). 
Suppose $A$ is a $\mathcal P$-algebra: this is 
a linear space equipped with bilinear operations 
$\mu: A\times A\to A$ for every $\mu \in V$.
These operations may be extended 
in a coefficient-wise way to the pairs of 
formal distributions over $A$ depending in different 
variables.

\begin{definition}
A pair of formal distributions 
$a(z),b(z)\in A[[z,z^{-1}]]$ is said to be {\em local}
if for every $\mu \in V$ there exists $N\in \mathbb Z_+$
such that 
\[
(w-z)^N \mu(a(w),b(z)) = 0.
\]
\end{definition}

In terms of coefficients, 
if $a(z)=\sum\limits_{n\in \mathbb Z} a(n)z^{-n-1}$,
$b(z)=\sum\limits_{m\in \mathbb Z} b(m)z^{-m-1}$,
$a(n),b(m)\in A$,
then 
the pair $a(z),b(z)$ is local if and only if
\begin{equation}\label{eq:LocalityCoeff}
\sum\limits_{s\ge 0} (-1)^s \binom{N}{s} \mu(a(n-s),b(m+s)) = 0, \quad N\gg 0,
\end{equation}
for all $n,m\in \mathbb Z$.
It is enough to check the locality condition for 
$\mu = e_i$, where $(e_i)_{i\in I}$ is a linear basis of $V$.

The locality condition naturally emerges in the theory 
of conformal algebras. Given a binary quadratic operad 
$\mathcal P$, suppose $M$ is a $\mathcal P$-conformal algebra, i.e., there is a morphism $\mathcal P\to \Cend_M$. If $(e_i)_{i\in I}$ is a basis of $V=\mathcal P(2)$ then denote by $*_i$ the image of $e_i$ in $\Cend_M(2)$:
\[
*_i: M\otimes M \to H^{\otimes 2}\otimes _H M, 
\]
where $H=\Bbbk [\partial ]$. Since $H$ is a Hopf algebra, one may uniquely rewrite $(a*_i b)$, $a,b\in M$, 
in the following form:
\[
(a*_i b) = \sum\limits_{n\ge 0} \dfrac{(-\partial)^n}{n!}\otimes c_{i(n)}, \quad c_{i(n)}\in M,
\]
for each $i\in I$ the sum is finite (see \cite{BDK2001}).
The elements $c_{i(n)}$ may be considered as the results 
of bilinear products $(a\ooc{i}{n} b)$ on $M$.

For every $\mathcal P$-conformal algebra $M$ defined as above there exists an ``ordinary'' $\mathcal P$-algebra $\mathcal A(M)$
with operations $\circ_i$, $i\in I$. 
Namely, the space 
$\mathcal A(M)$ is presented by $\Bbbk [t,t^{-1}]\otimes _H M$, where $\partial $ acts on $\Bbbk [t,t^{-1}]$
as $-d/dt$. 
Denote $t^n\otimes _H a = a(n)$ for $a\in M$,
$n\in \mathbb Z$. 
The products $\circ_i$ on $\mathcal A(M)$
are given by the well-defined formula
\[
a(n)\circ_i b(m) = \sum\limits_{s\ge 0} 
\binom{n}{s} (a\ooc{i}{s} b)(n+m-s), 
\quad a,b\in M,\ n,m\in \mathbb Z.
\]
These operations satisfy the defining identities 
of the variety of $\mathcal P$-algebras
(see \cite{Roit1999, Kol2006CA}). 

There exists an injective 
linear map 
\[
\iota : M \to \mathcal A(M)[[z,z^{-1}]], 
\quad 
a\mapsto \sum\limits_{n\in \mathbb Z} a(n)z^{-n-1}, 
\]
such that 
$\iota(\partial a) = \dfrac{d \iota(a)}{dz}$
for all $a\in M$ 
and 
\[
\iota(a\ooc{i}{n} b) = \iota(a)\ooc{i}{n} \iota(b),
\]
where the operations on the right-hand side are given 
by \eqref{eq:n-prod}.

For every $a,b\in M$, the formal distributions 
$\iota(a),\iota(b)$ form a local pair. 
Conversely, one may try to construct a $\mathcal P$-conformal algebra {\em generated} by a family of 
pairwise local formal distributions 
over an ordinary $\mathcal P$-algebra $A$. 
In order to get a success in this way, one needs the general statement similar to the Dong Lemma. 

\begin{definition}\label{defn:DongOperad}
Let us say that a binary quadratic 
operad $\mathcal P=\mathcal P(V,R)$ satisfies 
the Dong Property if for every $\mathcal P$-algebra 
$A$ and for every pairwise local formal distributions
$a(z),b(z),c(z) \in A[[z,z^{-1}]]$
the distributions $(a\ooc{i}{n} b)(z)$ and $c(z)$ 
form a local pair for every $i\in I$, $n\in \mathbb Z_+$.
\end{definition}

Here, as above, we have fixed a linear basis $(e_i)_{i\in I}$ of $V$. Obviously, the locality property does not depend on the choice of a basis, 
so the Dong Property is immanent to the operad $\mathcal P$. Alternatively, we will say $\mathcal P$ is a {\em Dong operad}.

The main purpose of this paper is to set up a criterion for a binary quadratic operad $\mathcal P$ to satisfy the Dong Property. This criterion will be stated in terms of the Koszul dual operad $\mathcal P^!$. 
Recall that $\mathcal P^!=\mathcal P(V^\vee, R^\perp)$, where $R^\perp \subseteq (\mathcal F_V(3))^\vee
\cong \mathcal F_{V^\vee}(3)$ is the subspace orthogonal to~$R$.

It turns out that the Dong Property 
for $\mathcal P$ is determined by the intersection 
of the subspaces 
$R^\perp$ and $V^\vee \vecotimes V^\vee$ in 
$\mathcal F_{V^\vee}(3)$.
Namely, we will be interested in the condition 
\begin{equation}\label{eq:CriterionDong}
(V^\vee \vecotimes V^\vee)\cap R^\perp = 0
\end{equation}
which is easy to check provided that we know 
a basis of $\mathcal P^!(3)$.
Fix a linear basis $(e_i)_{i\in I}$ of $V$
and choose the dual basis $(e_i^\vee )_{i\in I}$
of the transpose dual space $V^\vee$.

Suppose $A$ is a $\mathcal P$-algebra with operations 
$e_i(a,b) = a\circ_i b$, $i\in I$, $a,b\in A$.
Similarly, every algebra $B$ from the class of
$\mathcal P^!$-algebras is equipped with operations 
$e_i^\vee (x,y) = x\circ^i y$, $i\in I$, $x,y\in B$.

The condition \eqref{eq:CriterionDong} is equivalent to the following one: 
in the free algebra $\mathcal P^!\langle p,q,t\rangle $ generated by a 3-element set $\{p,q,t\}$,
the elements
\begin{equation}\label{eq:NSymmElements}
(p\circ^i q)\circ^j t, \quad i,j\in I,
\end{equation}
are linearly independent (they correspond to 
$e_j^\vee \vecotimes e_i^\vee \in V^\vee \vecotimes V^\vee $). 

In particular, if $\dim V=1$, $e_1(a,b) = ab$,
as in the commutative or anti-commutative cases,
then the family \eqref{eq:NSymmElements} contains 
only one term $(pq)t$. Hence, 
\eqref{eq:CriterionDong} is equivalent to 
$(pq)t\ne 0$ in the free $\mathcal P^!$-algebra.

If $V\cong S_2$, as in those cases like associative, 
alternative, Novikov, or Poisson algebras, then 
we may choose $e_2 = (12)e_1$, so that 
$e_1(a,b)=ab$, $e_2(a,b)=ba$. In this case 
the family \eqref{eq:NSymmElements} contains four terms:
\[
(pq)t,\ (qp)t, \ t(pq),\ t(qp).
\]
Therefore, the condition \eqref{eq:CriterionDong}
means that these four monomials are linearly independent 
in the free $\mathcal P^!$-algebra.

\begin{theorem}\label{thm:Sufficient}
Assume that for a binary quadratic operad
$\mathcal P=\mathcal P(V,R)$ 
we have $(V^\vee \vecotimes V^\vee)\cap R^\perp = 0$.
Then $\mathcal P$ is a Dong operad.
\end{theorem}

\begin{proof}
Suppose $A$ is a $\mathcal P$-algebra with operations 
$(\cdot \circ_i \cdot)$ as above, and let $B$ be the free 
$\mathcal P^!$-algebra generated by $\{p,q,t\}$
with dual operations $(\cdot \circ^i \cdot )$.

Assume $a(z),b(z),c(z) \in A[[z,z^{-1}]]$
are three pairwise local formal distributions over $A$,
\[
x(z) = \sum\limits_{n\in \mathbb Z} x(n)z^{-n-1},
\quad x(n)\in A,\ x\in \{ a,b,c\}.
\]
Denote
\[
\bar a(z) = p\otimes a(z) = \sum\limits_{n\in \mathbb Z} p\otimes a(n)
z^{-n-1} 
\] 
is a formal distribution over $B\otimes A$. 
Similarly, define $\bar b(z) = q\otimes b(z)$
and $\bar c(z)=t\otimes c(z)$.

Recall that $B\otimes A$ is a Lie algebra relative to the operation \eqref{eq:DualProdLie}.

\begin{lemma}
The formal distributions 
$\bar a(z),\bar b(z), \bar c(z) \in in (B\otimes A)[[z,z^{-1}]]$ are pairwise local.
\end{lemma}

\begin{proof}
Consider the pair $\bar a(z),\bar b(z)$.
By definition, there exists $N\ge 0$ such that 
\[
\sum\limits_{s\ge 0} (-1)^s \binom{N}{s}
 a(n-s)\circ_i b(m+s) = 0 \in A
\]
for all $n,m\in \mathbb Z$, $i\in I$
(recall that $|I|=\dim V<\infty $).
The Lie bracket on $B\otimes A$
is given by \eqref{eq:DualProdLie}, 
so 
\[
[p\otimes a(n), q\otimes b(m)]
=
\sum\limits_{i\in I}
(p\circ^i q) \otimes ( a(n)\circ_i b(m) )
\]
for $n,m\in \mathbb Z$.
Hence, for the same $N$ we have
\begin{multline*}
\sum\limits_{s\ge 0} (-1)^s \binom{N}{s}
[p\otimes a(n-s), q\otimes b(m+s)]
=
\sum\limits_{s\ge 0} (-1)^s \binom{N}{s}
\sum\limits_{i\in I}
(p\circ^i q) \otimes ( a(n-s)\circ_i b(m+s) ) \\
=
\sum\limits_{i\in I}
(p\circ^i q) \otimes 
\sum\limits_{s\ge 0} (-1)^s \binom{N}{s}
 a(n-s)\circ_i b(m+s) 
 = 0,
\end{multline*}
as desired.
\end{proof}

Hence, the Dong Lemma for Lie algebras
can be applied to the distributions $\bar a(z),\bar b(z), \bar c(z) \in (B\otimes A)[[z,z^{-1}]]$.
Namely,  for every $k\in \mathbb Z_+$
the formal distributions 
$[\bar a\oo{k}\bar b](z)$ and $\bar c(z)$
are local.

The coefficient of $[\bar a\oo{k} \bar b](z)$ 
at $z^{-n-1}$
is given by 
\eqref{eq:n-prod}, that is,
\[
[\bar a\oo{k}\bar b](n) = 
\sum\limits_{t\ge 0} (-1)^t\binom{k}{t}
[p\otimes a(k-t), q\otimes b(n+t)].
\]

Hence, 
there exists $N\in \mathbb Z_+$
such that 
\begin{equation}\label{eq:LieDong-coeff}
\sum\limits_{t,s\ge 0}
(-1)^{s+t}\binom{k}{t}\binom{N}{s}
[ [p\otimes a(k-t), q\otimes b(n-s+t)], t\otimes c(m+s) ] = 0
\end{equation}
for every $n,m\in \mathbb Z$.
Expand \eqref{eq:LieDong-coeff} 
by definition \eqref{eq:DualProdLie}:
\begin{multline*}
\sum\limits_{t,s\ge 0}
\sum\limits_{i\in I}
(-1)^{t+s} \binom{k}{t} \binom{N}{s}
[(p\circ ^i q)\otimes 
(a(k-t)\circ_i b(n-s+t)),
t\otimes c(m+s) ]
\\
=
\sum\limits_{i,j\in I}
\sum\limits_{t,s\ge 0}
(-1)^{t+s} \binom{k}{t} \binom{N}{s}
((p\circ ^i q)\circ^j t)\otimes 
(a(k-t)\circ_i b(n-s+t)) \circ_j c(m+s) 
=0.
\end{multline*}
Since $(p\circ^i q)\circ^j t$, $i,j\in I$, 
are linearly 
independent, the corresponding coefficients 
are zero in $A$:
\[
\sum\limits_{t,s\ge 0}
(-1)^{t+s} \binom{k}{t} \binom{N}{s}
(a(k-t)\circ_i b(n-s+t)) \circ_j c(m+s) 
=0
\]
for all $i,j\in I$. The latter means 
the distributions 
$(a\ooc{i}{k} b)(z)$ and $c(z)$ 
form a local pair. 
\end{proof}

\section{The Dong Property: necessary condition}

According to Theorem~\ref{thm:Sufficient},
the condition
\eqref{eq:CriterionDong}
is sufficient for a binary quadratic operad 
to satisfy the Dong Property in the sense of 
Definition~\ref{defn:DongOperad}.
The purpose of this section is to prove that 
\eqref{eq:CriterionDong}
is also a necessary condition. 

\begin{theorem}\label{thm:Necessary}
Suppose $\mathcal P=\mathcal P(V,R)$ 
is a binary quadratic operad.
If $\mathcal P$ satisfies the Dong Property then 
the canonical epimorphism 
\[
\mathcal F_{V^\vee}(3) \to \mathcal F_V(3)^\vee/R^\perp = \mathcal P^!(3)
\]
acts injectively on
$V^\vee \vecotimes V^\vee $.
\end{theorem}

\begin{proof}
Let us describe the details of the proof 
in a particular but important case when 
$V\simeq \Bbbk S_2$, i.e., when 
$\mathcal P$ is generated by single non-symmetric 
operation $(x,y)\mapsto xy$
(then so is the dual operad $\mathcal P^!$). 

Assume the conclusion is wrong, i.e., there exists 
nonzero $h\in V^\vee \vecotimes V^\vee $
such that the image of $h$ is zero in $\mathcal P^!(3)$. 
The latter means, by definition, that $h \in R^\perp$, 
where $R$ is the $S_3$-submodule of defining relations of 
$\mathcal P$.
In the generic form, $h$ is a linear combination 
of $(x_1x_2)x_3$, $(x_2x_1)x_3$, $x_3(x_1x_2)$, and $x_3(x_2x_1)$. There are two different cases up to permutation $(12)$.

{\sc Case 1.} A coefficient at $(x_1x_2)x_3$ in $h$ is nonzero, 
i.e., 
\begin{equation}\label{eq:H-polynomial}
h =  (x_1x_2)x_3 + \alpha (x_2x_1)x_3 + \beta x_3(x_1x_2) +\gamma x_3(x_2x_1),
\quad \alpha,\beta,\gamma \in \Bbbk .
\end{equation}

{\sc Case 2.} A coefficient at $x_3(x_1x_2)$ in $h$ 
is nonzero.

Let us consider Case~1 in details: we will construct 
a $\mathcal P$-algebra $A$ and three formal distributions 
$a,b,c\in A[[z,z^{-1}]]$ such that 
$a,b,c$ are pairwise local but the pair $(a\oo 0 b, c)$
is not a local one.
(For Case 2, a contradiction can be obtained 
similarly, just the non-local pair is $(a, b\oo 0 c)$.)

Denote by $X$ the set of all symbols $x(n)$, where 
$n\in \mathbb Z$, $x\in \{a,b,c\}$. 
Consider the free $\mathcal P$-algebra 
$\mathcal P\langle X\rangle $
generated by 
$X$
and let $I$ stand for the ideal of 
$\mathcal P\langle X\rangle $ 
generated by 
the relations 
\begin{equation}\label{eq:Rel_Loc1}
x(n)y(m)-x(n-1)y(m+1),\quad x,y\in \{a,b,c\},\ n,m\in \mathbb Z.
\end{equation}
Then the $\mathcal P$-algebra 
$A=\mathcal P\langle X\rangle /I$ has the following property:
the formal distributions 
\[
x(z) = \sum\limits_{n\in \mathbb Z} x(n)z^{-n-1}\in A[[z,z^{-1}]], 
\quad x\in \{a,b,c\},
\]
are pairwise local with $N=1$.

Since the operad $\mathcal P$ has the Dong Property, 
the formal distributions $(a\oo 0 b)(z)$ and $c(z)$
also form a local pair, i.e., there exists 
$N\in \mathbb Z_+$ such that 
\[
f=\sum\limits_{s\ge 0} (-1)^s \binom{N}{s}
(a(0)b(n-s))c(m+s) \in I \triangleleft \mathcal P\langle X\rangle 
\]
for all $n,m\in \mathbb Z$.
The latter means that $f$ may be presented as a linear combination of left and right multiples of the relations
\eqref{eq:Rel_Loc1}, i.e.,
\begin{multline}\label{eq:Generic_Loc}
f = \sum\limits_{x,y,z\in \{a,b,c\}}
\sum\limits_{i+j+k=m+n} \gamma_{x,y,z}^{i,j,k}
(x(i)y(j)-x(i-1)y(j+1))z(k) \\
+
\sum\limits_{x,y,z\in \{a,b,c\}}
\sum\limits_{i+j+k=m+n} \beta_{x,y,z}^{i,j,k}
z(k)(x(i)y(j)-x(i-1)y(j+1))
\end{multline}
in the free $\mathcal P$-algebra.
Note that the sums are finite, so almost all 
$\gamma_{x,y,z}^{i,j,k}, \beta_{x,y,z}^{i,j,k} \in \Bbbk $
are zero.

Let us re-arrange the terms in the right-hand side 
of \eqref{eq:Generic_Loc} to get
\begin{multline}\label{eq:Generic_Loc-2}
-\sum\limits_{s\ge 0}
(-1)^s \binom{N}{s} (a(0)b(n-s))c(m+s) 
+ \sum\limits_{x,y,z\in \{a,b,c\}}
\sum\limits_{i+j+k=m+n} 
(\gamma_{x,y,z}^{i,j,k} -\gamma_{x,y,z}^{i+1,j-1,k})
(x(i)y(j))z(k) \\
+
\sum\limits_{x,y,z\in \{a,b,c\}}
\sum\limits_{i+j+k=m+n} 
(\beta_{x,y,z}^{i,j,k} - \beta_{x,y,z}^{i+1,j-1,k})
z(k)(x(i)y(j)) =0 \in \mathcal P\langle X\rangle .
\end{multline}

Given a triple of integers $i,j,k$, denote by 
$f^{i,j,k}$ the sum of all those summands 
of \eqref{eq:Generic_Loc-2}
that contain the variables $a(i),b(j),c(k)$. 
All polynomials of the form $f^{i,n-i,m}$, 
$i\in \mathbb Z$, are of degree~3 and they are
equal to zero in the free 
$\mathcal P$-algebra. 
Hence, each  $f^{i,n-i,m}$ is obtained from some element of $R$ 
by replacing $x_1,x_2,x_3$ by the respective variables 
$a(i),b(n-i),c(m)\in X$, i.e.,  
\[
F^{i,n-i,m}|_{x_1=a(i),x_2=b(n-i),x_3=c(m)} .
\]
In particular,
\begin{multline}\label{eq:F0nm}
f^{0,n,m} 
= -(a(0)b(n))c(m) 
+ (\gamma_{a,b,c}^{0,n,m} - \gamma_{a,b,c}^{1,n-1,m})
  (a(0)b(n))c(m)
+ (\gamma_{b,a,c}^{n,0,m} - \gamma_{b,a,c}^{n+1,-1,m})
  (b(n)a(0))c(m)
  \\
+ (\gamma_{a,c,b}^{0,m,n} - \gamma_{a,c,b}^{1,m-1,n})
  (a(0)c(m))b(n)
+ (\gamma_{c,a,b}^{m,0,n} - \gamma_{c,a,b}^{m+1,-1,n})
  (c(m)a(0))b(n)
  \\
+ (\gamma_{b,c,a}^{n,m,0} - \gamma_{b,c,a}^{n+1,m-1,0})
  (b(n)c(m))a(0)
+ (\gamma_{c,b,a}^{m,n,0} - \gamma_{c,b,a}^{m+1,n-1,0})
  (c(m)b(n))a(0)
\\
+ (\beta_{a,b,c}^{0,n,m} - \beta_{a,b,c}^{1,n-1,m})
  c(m)(a(0)b(n))
+ (\beta_{b,a,c}^{n,0,m} - \beta_{b,a,c}^{n+1,-1,m})
  c(m)(b(n)a(0))
\\
+ (\beta_{a,c,b}^{0,m,n} - \beta_{a,c,b}^{1,m-1,n})
  b(n)(a(0)c(m))
+ (\beta_{c,a,b}^{m,0,n} - \beta_{c,a,b}^{m+1,-1,n})
  b(n)(c(m)a(0))
\\
+ (\beta_{b,c,a}^{n,m,0} - \beta_{b,c,a}^{n+1,m-1,0})
  a(0)(b(n)c(m))
+ (\beta_{c,b,a}^{m,n,0} - \beta_{c,b,a}^{m+1,n-1,0})
  a(0)(c(m)b(n)).
\end{multline}

By definition of a dual operad, 
$\langle h,F^{0,n,m} \rangle =0$, 
where $h$ is the element 
of $\mathcal F_{V^\vee}(3)$ represented by 
\eqref{eq:H-polynomial}. 
The latter pairing is 
easy to compute in terms of the coefficients from 
\eqref{eq:F0nm}:
\begin{multline*}
0 = \langle h,F^{0,n,m} \rangle \\
= 
-1 
+ (\gamma_{a,b,c}^{0,n,m} - \gamma_{a,b,c}^{1,n-1,m})
-\alpha (\gamma_{b,a,c}^{n,0,m} - \gamma_{b,a,c}^{n+1,-1,m})
-\beta (\beta_{a,b,c}^{0,n,m} - \beta_{a,b,c}^{1,n-1,m})
+ \gamma (\beta_{b,a,c}^{n,0,m} - \beta_{b,a,c}^{n+1,-1,m}),
\end{multline*}
i.e., 
\begin{equation}\label{eq:Zero0-orth}
(\gamma_{a,b,c}^{0,n,m} - \gamma_{a,b,c}^{1,n-1,m})
-\alpha (\gamma_{b,a,c}^{n,0,m} - \gamma_{b,a,c}^{n+1,-1,m})
-\beta (\beta_{a,b,c}^{0,n,m} - \beta_{a,b,c}^{1,n-1,m})
+ \gamma (\beta_{b,a,c}^{n,0,m} - \beta_{b,a,c}^{n+1,-1,m})    =1.
\end{equation}
For $F^{1,n-1,m}$ we also have 
$\langle h, F^{1,n-1,m}\rangle =0$, 
but note that $a(1)$ does not appear in the 
locality relation~$f$.  Hence, 
\begin{equation}\label{eq:Zero1-orth}
(\gamma_{a,b,c}^{1,n-1,m} - \gamma_{a,b,c}^{2,n-2,m})
-\alpha (\gamma_{b,a,c}^{n-1,1,m} - \gamma_{b,a,c}^{n,0,m})
-\beta (\beta_{a,b,c}^{1,n-1,m} - \beta_{a,b,c}^{2,n-2,m})
+ \gamma (\beta_{b,a,c}^{n-1,1,m} - \beta_{b,a,c}^{n,0,m})    =0.
\end{equation}
Now, let us add \eqref{eq:Zero0-orth} and \eqref{eq:Zero1-orth} to get
\begin{equation}\label{eq:Zero01-orth}
(\gamma_{a,b,c}^{0,n,m} - \gamma_{a,b,c}^{2,n-2,m})
-\alpha (\gamma_{b,a,c}^{n-1,1,m} - \gamma_{b,a,c}^{n+1,-1,m})
-\beta (\beta_{a,b,c}^{0,n,m} - \beta_{a,b,c}^{2,n-2,m})
+ \gamma (\beta_{b,a,c}^{n-1,1,m} - \beta_{b,a,c}^{n+1,-1,m})    =1.
\end{equation}
Similarly, for every $i=2,3,\dots $ it follows 
from $\langle h, F^{i,n-i,m}\rangle =0$ that
\begin{multline}\label{eq:Zero2-orth}
(\gamma_{a,b,c}^{i,n-i,m} - \gamma_{a,b,c}^{i+1,n-i-1,m})
-\alpha (\gamma_{b,a,c}^{n-i,i,m} - \gamma_{b,a,c}^{n-i+1,i-1,m}) 
\\
-\beta (\beta_{a,b,c}^{i,n-i,m} - \beta_{a,b,c}^{i+1,n-i-1,m})
+ \gamma (\beta_{b,a,c}^{n-i,i,m} - \beta_{b,a,c}^{n-i+1,i-1,m})    =0.
\end{multline}
Let us subsequently add \eqref{eq:Zero2-orth} for 
$i=2,3,\dots, k$, $k\in \mathbb Z$, to 
\eqref{eq:Zero01-orth} to obtain
\begin{multline}\label{eq:Zero0k-orth}
(\gamma_{a,b,c}^{0,n,m} - \gamma_{a,b,c}^{k,n-k,m})
-\alpha (\gamma_{b,a,c}^{n-k+1,k-1,m} - \gamma_{b,a,c}^{n+1,-1,m})
\\
-\beta (\beta_{a,b,c}^{0,n,m} - \beta_{a,b,c}^{k,n-k,m})
+ \gamma (\beta_{b,a,c}^{n-k+1,k-1,m} - \beta_{b,a,c}^{n+1,-1,m})    =1.
\end{multline}
Since almost all coefficients in the decomposition 
\eqref{eq:Generic_Loc} are zero, for $k\gg 0$ we have
\begin{equation}\label{eq:PositiveReduce}
\gamma_{a,b,c}^{0,n,m}
+\alpha \gamma_{b,a,c}^{n+1,-1,m}
-\beta \beta_{a,b,c}^{0,n,m}
- \gamma \beta_{b,a,c}^{n+1,-1,m}    =1.
\end{equation}

Now, deduce similar relations that follow from 
$\langle h, F^{-i,n+i,m}\rangle =0 $, $i=1,2,\dots $:
for example, if $i=1$ then 
\begin{multline}\label{eq:ZeroMinus-orth}
(\gamma_{a,b,c}^{-1,n+1,m} - \gamma_{a,b,c}^{0,n,m})
-\alpha (\gamma_{b,a,c}^{n+1,-1,m} - \gamma_{b,a,c}^{n+2,-2,m}) 
\\
-\beta (\beta_{a,b,c}^{-1,n+1,m} - \beta_{a,b,c}^{
0,n,m})
+ \gamma (\beta_{b,a,c}^{n+1,-1,m} 
- \beta_{b,a,c}^{n+2,-2,m})    =0.
\end{multline}
If we add \eqref{eq:ZeroMinus-orth} to \eqref{eq:PositiveReduce} then the relation obtained 
is
\[
\gamma_{a,b,c}^{-1,n+1,m}
+\alpha \gamma_{b,a,c}^{n+2,-2,m}
-\beta \beta_{a,b,c}^{-1,n+1,m}
- \gamma \beta_{b,a,c}^{n+2,-2,m}    =1.
\]
Proceed in a similar way to get 
\begin{equation}\label{eq:NegativeReduce}
\gamma_{a,b,c}^{-k,n+k,m}
+\alpha \gamma_{b,a,c}^{n+k+1,-k-1,m}
-\beta \beta_{a,b,c}^{-k,n+k,m}
- \gamma \beta_{b,a,c}^{n+k+1,-k-1,m}    =1.
\end{equation}
for every $k\ge 1$. When $k$ is large enough,
we obtain $0=1$, a contradiction.
\end{proof}

\section{Applications to Manin products}

By Theorem \ref{thm:Sufficient}, the condition 
\eqref{eq:CriterionDong}
is sufficient for a binary quadratic operad 
$\mathcal P=\mathcal P(V,R)$ to satisfy the Dong Property. 
Apart from well-known examples of such operads 
($\mathcal P = \As,\Com,\Lie$), this sufficient 
condition holds for many other classes of algebras.

\begin{example}
Let $\mathcal P=\Pois$ be the operad of Poisson 
algebras (see Example \ref{exmp:Pois}). 
Then $\Pois^!$ is isomorphic to $\Pois $, so this is a Dong operad.
\end{example}

Indeed, the family \eqref{eq:NSymmElements} for Poisson algebras consists of four elements 
\[
\{\{p,q\},t\}, \ \{p,q\}t,\ \{pq,t\}, \ pqt.
\]
They are linearly independent in the free Poisson algebra.

\begin{example}
Let $\mathcal P=\Nov$ be the operad of Novikov algebras 
from Example \ref{exmp:Nov}.
Then the Koszul dual operad $\Nov^!$ describes 
the class of right-symmetric and left commutative algebras, so the condition of Theorem~\ref{thm:Sufficient} holds, and
$\Nov $ is a Dong operad.
\end{example}

Indeed, the family \eqref{eq:NSymmElements} may be presented by elements from the free commutative differential  algebra \cite{DzhLofwall2002}:
\[
(x\circ y)\circ z = x''yz+x'y'z,
\ 
(y\circ x)\circ z = xy''z+x'y'z,
\ 
z\circ (x\circ y) = x'yz',
\ 
z\circ (y\circ x) = xy'z',
\]
they are linearly independent. Therefore, the class of Novikov algebras satisfies the Dong Property (in a different way it was proved in \cite{BokKol2024}).

\begin{example}
Let $\mathcal P=\NP$ be the operad of Novikov--Poisson 
algebras from Example~\ref{exmp:NovPois}.
It follows from the description of $\NP^!$ (see Example~\ref{exmp:NP-dual}) 
that the elements \eqref{eq:NSymmElements} 
\[
[[p,q],t],\ [p,q]\circ t, \ t\circ [p,q],
\ [p\circ q, t],\ [q\circ p, t], 
\ (p\circ q)\circ t, \ (q\circ p)\circ t,
\ t\circ (p\circ q), \ t\circ (q\circ p)
\]
are linearly independent in the free $\NP^!$-algebra.
Hence, $\NP$ is a Dong operad.
\end{example}

\begin{example}\label{exmp:Alt}
Let $\mathcal P = \mathrm{Alt}$ be the operad of alternative algebras, i.e., there is one binary product 
satisfying the identities
\[
\begin{gathered}
(x_1x_2)x_3 - x_1(x_2x_3) = -(x_2x_1)x_3 + x_2(x_1x_3), \\
(x_1x_2)x_3 - x_1(x_2x_3) = -(x_1x_3)x_2 +x_1(x_3x_2). 
\end{gathered}
\]
\end{example}

Then $\mathcal P^!$ is the operad of associative algebras satisfying $x^3=0$. The latter identity 
in the linearized form is
\[
\sum\limits_{\sigma \in S_3}x_{\sigma(1)}x_{\sigma(2)}x_{\sigma(3)},
\]
so the terms $x_1x_2x_3, x_2x_1x_3, x_3x_1x_2, x_3x_2x_1$
are linearly independent in $\mathcal P^!(3)$.
Therefore, the Dong Property holds for alternative algebras.

However, $\NP^!$ is not a Dong operad, so the more 
general class of $\GD $-algebras also does not satisfy the Dong Property. Note that the dual operad $\GD^!$
is Dong. 

Since the paper \cite{LodayPirash199?}, various classes 
of di- and tri-algebras have been widely studied
(see \cite{BGG_HK_Proceed} and references therein).
For every variety $\Var $ of algebras one may construct 
the classes $\di\Var$ and $\tri\Var $ that are 
defined by the operads $\Perm\otimes \Var$ and 
$\ComTriAs\otimes \Var $. For a binary quadratic operad
$\mathcal P $,  the Hadamard product $\Perm\otimes \mathcal P$ 
coincides with the so-called Manin white product
\cite{GinKapr1994}
denoted
$\Perm\circ \mathcal P$.
The same holds for the product 
with $\ComTriAs$. In particular, $\di\Lie $ is the operad of Leibniz algebras.

The natural question is to decide whether 
$\di\mathcal P $
or $\tri\mathcal P $ 
is a Dong operad provided that so is $\mathcal P$.
Since the criterion \eqref{eq:CriterionDong} concerns 
the Koszul dual operad, we have to study 
\[
(\di\mathcal P)^! = \pre\Lie\bullet \mathcal P^!, \quad 
(\tri\mathcal P)^! = \post\Lie\bullet \mathcal P^!, 
\]
where $\bullet $ stands for the Manin black product of \cite{GinKapr1994}, 
$\pre\Lie $ and $\post\Lie $ are the operads of pre-Lie and post-Lie algebras, respectively. 

\begin{corollary}
Let $\mathcal P$ be a Dong operad. Then 
the operads $\di\mathcal P$, $\tri\mathcal P$ governing the classes of replicated 
algebras are also Dong.
\end{corollary}

\begin{proof}
Given a multilinear variety $\Var $, the identities 
of $\pre\Var$ or $\post\Var$ can be found by means 
of an algorithmic procedure described in \cite{BBGN13}
(see also \cite{GK13, GK14}).

If $\Var$ is the class of $\mathcal P^!$-algebras, where 
$\mathcal P=\mathcal P(V,R)$ is a binary quadratic operad, 
then the algebras in $\Var $ are linear spaces equipped with a family of operations 
$e_i^\vee = x_1\circ ^i x_2$, $(e_i)_{i\in I}$ is a basis of $V$. The algebras from $\pre\Var $
have duplicated family of operations 
$\mu_1\otimes e_i^\vee = x_1\succ^i x_2$,
$\mu_2\otimes e_i^\vee = x_1\prec^i x_2$, $i\in I$,
where $\mu_1,\mu_2$ form a basis 
of operations for $\pre\Lie$, 
$(12)\mu_1=\mu_2$. 
Similarly, the algebras from 
$\post\Var $ have triple family of operations:
$\mu_1\otimes e_i^\vee = x_1\succ^i x_2$,
$\mu_2\otimes e_i^\vee = x_1\prec^i x_2$
as above, 
and 
$\nu\otimes e_i^\vee = x_1\perp^i x_2$. 
Here $\mu_1,\mu_2,\nu $ form a basis of the operation
space for $\post\Lie$, $(12)\mu_1=\mu_2$, $(12)\nu=-\nu$.

Let us consider the tri-case: recall how to find 
the identities of $\post\Var$ \cite{GK14}. 
For brevity, denote 
\[
a*^j b = a\succ^j b+a\prec^j b  +a\perp^j b.
\]
For each defining identity $f(x_1,x_2,x_3)$ of $\Var $
(these identities form the space $R^\perp $) consider the {\em splitting} $f_1,\dots,f_7$ constructed as follows. Suppose $f$ is a linear combination 
of monomials in $x_1,x_2,x_3$ with respect 
to the operations $\circ^i$, $i\in I$.
For each monomial 
$u =(x_{k_1}\circ ^j x_{k_2})\circ^i x_{k_3}$, 
$\{k_1,k_2,k_3\}=\{1,2,3\}$, $i,j\in I$,
consider 
\[
\begin{gathered}
 u_{k_1} = (x_{k_1}\prec ^j x_{k_2})\prec^i x_{k_3}, 
 \quad 
 u_{k_2} = (x_{k_1}\succ ^j x_{k_2})\prec^i x_{k_3}, 
 \quad
 u_{k_3} = (x_{k_1} *^j x_{k_2})\succ^i x_{k_3}, 
 \\
 u_{k_1,k_2} = (x_{k_1}\perp^j x_{k_2})\prec^i x_{k_3}, 
 \quad
 u_{k_1,k_3} = (x_{k_1}\prec ^j x_{k_2})\perp ^i x_{k_3}, 
 \quad
 u_{k_2,k_3} = (x_{k_1}\succ ^j x_{k_2})\perp^i x_{k_3}, 
 \\
 u_{k_1,k_2,k_3} = (x_{k_1}\perp^j x_{k_2})\perp^i x_{k_3}.
\end{gathered}
\]
We do not need to consider the right-justified 
monomials
$u =x_{k_1}\circ ^i (x_{k_2} \circ^j x_{k_3})$
separately since we assume the family $e_i^\vee$
already contains permutations $(12)e_i^\vee$.

In order to find the defining identities of $\post\Var$
coming from an identity $f$ of $\Var$
one has to replace each monomial $u$ in $f$ with the respective $u_{M}$, $M$ is a non-empty subset of
$\{1,2,3\}$. The obtained relations $f_M$ form the desired splitting.

Denote by $R_M^\perp $ the linear span of all $f_M$, 
$f\in R^\perp$, $\emptyset M\subseteq \{1,2,3\}$. 
This is an $S_3$-subspace of $\mathcal F_{W}(3)$, 
where $W=\post\Lie(2)\otimes V^\vee $.

For every $f\in \mathcal F_{V^\vee }(3)$ we may write 
\[
f = f^{(0)} + (13)f^{(1)} + (23)f^{(2)}, 
\]
where $f^{(i)}\in V^\vee\vecotimes V^\vee $.
Note that the splitting procedure described above 
does not affect the order of variables $x_1,x_2,x_3$.
Hence, if the elements $(13)f^{(1)} + (23)f^{(2)}$, 
$f\in R^\perp $, were linearly independent in 
$ \mathcal F_{V^\vee }(3)$ then so are their splittings
as summands of 
\[
f_M = f_M^{(0)} + (13)f_M^{(1)} + (23)f_M^{(2)}.
\]
Therefore, every nonzero polynomial $g$ from 
$R_M^\perp $ contains a nonzero summand $g^{(1)}$ or $g^{(2)}$ provided that \eqref{eq:CriterionDong} holds
for the initial operad~$\mathcal P$. 
\end{proof}

The next natural question is to study how the Dong Property relates with the splitting procedure.
It was shown in \cite{BokKol2024} that the operads
$\pre\As$ and $\pre\Lie$ are not Dong.
Now we may check it again using the criterion 
\eqref{eq:CriterionDong}.

\begin{example}\label{exmp:preAs}
Let $\mathcal P=\pre\As$ be the operad of pre-associative (dendriform) algebras. Then 
$\mathcal P^!=\di\As$ is the operad of di-associative
algebras (dialgebras).
\end{example}

Choose a basis of $\di\As(2)$ to be of 
$e_1 = x_1\vdash x_2$,
$e_2 =(12)e_1 =x_2\vdash x_1$,
$e_3 =x_1\dashv x_2$,
$e_4 =(12)e_3 =x_2\dashv x_1$.
If we assume $\pre\As$ satisfies the Dong Property 
then 
\[
e_i\vecotimes e_j, \quad i,j=1,\ldots,4,
\]
are linearly independent in $\di\As(3)$. 
However, we have 
\[
e_1\vecotimes e_1 - e_1\vecotimes e_3 = (x_1\vdash x_2)\vdash x_3 - (x_1\dashv x_2)\vdash x_3 =0 \in \di\As(3), 
\]
i.e., the criterion of Theorem~\ref{thm:Necessary}
does not hold. Hence, the operad $\pre\As$ does not satisfy 
the Dong Property.

Even the smaller variety of pre-commutative (Zinbiel)
algebras does not have the Dong Property since 
$\Zinb^! = \Leib$.
By similar reason, since $(x_1x_2)x_3-(x_2x_1)x_3=0\in \Perm(3)$, the operad $\pre\Lie=\Perm^!$
does not satisfy the Dong Property.

In general, we may deduce the following

\begin{corollary}
Let $\mathcal P$ be a binary quadratic 
operad. Then the dendriform splitting operad $\pre\mathcal P = \pre\Lie \bullet \mathcal P$
does not satisfy the Dong Property.
\end{corollary}

\begin{proof}
Indeed, $(\pre\mathcal P)^! = \di(\mathcal P^!) = \Perm\otimes \mathcal P$ \cite{BBGN13, GK13}.
Suppose 
$\mathcal P =\mathcal P(V,R)$,
and $0\ne \mu\in V^\vee $. 
Then for $x_1x_2\in \Perm(2)$ we have 
\begin{multline*}
(x_1x_2\otimes \mu)\vec\otimes (x_1x_2\otimes \mu )
-
(x_1x_2\otimes \mu)\vec\otimes (x_2x_1\otimes \mu )
\\
=
(x_1x_2)x_3\otimes (\mu\vec\otimes \mu) -
(x_2x_1)x_3\otimes (\mu\vec\otimes \mu) =0\in 
\Perm(3)\otimes \mathcal P^!(3) = \di(\mathcal P^!)(3).
\end{multline*}
Hence the necessary condition for the operad 
$\pre\mathcal P$ to satisfy the Dong Property 
does not hold.
\end{proof}

Obviously, the greater operads (tridendriform splitting) $\post\Var=\post\Lie\bullet \Var $ also 
do not satisfy the Dong Property.

However, if we are given two Dong operads $\mathcal P$
and $\mathcal Q$ then it turns out that their Manin black product $\mathcal P\bullet \mathcal Q$ 
is also a Dong operad (this statement is not true for the white product).

In order to prove that, we need an algorithmic way 
to compute the defining relations of $\mathcal P\bullet \mathcal Q$. It was actually mentioned in \cite[2.2.8]{GinKapr1994}, but in a very sketched form. In \cite{GK14} this method (with an independent proof based on Rota--Baxter operators) was applied to describe the dendriform splitting of an operad $\mathcal Q$, i.e., to find the defining relations of $\pre\Lie\bullet\mathcal Q$ or $\post\Lie\bullet\mathcal Q$. Let us present an easy way to find the Manin black product of binary quadratic operads in general.

Let $\mathcal P=\mathcal P(V,R)$ and 
$\mathcal Q=\mathcal P(W,S)$ 
be two binary quadratic operads,
$\dim V,\dim W<\infty $.
Assume $(e_i)_{i\in I}$ and $(f_j)_{j\in J}$
are linear bases of $V$ and $W$, respectively.
Choose the dual basis $(e_i^\vee)_{i\in I}$ of $V^\vee $.

As above, express the operations generating 
$\mathcal P^!$ and $\mathcal Q$ 
via monomials: 
\[
e^\vee_i = x_1\circ^i x_2,\quad f_j = x_1 *_j x_2.
\]
Suppose $A$ is an algebra with the space of binary operations presented by 
$(1_-\otimes e_i\otimes f_j)_{(i,j)\in I\times J}$,
\[
1_-\otimes e_i\otimes f_j = x_1 *_{i,j} x_2, \quad i\in I,\ j\in J.
\]
Then for every $\mathcal P^!$-algebra $B$
the space $B\otimes A$ may be considered 
as an algebra with operations $(\cdot *_j \cdot)$,
$j\in J$,
defined as follows:
\begin{equation}\label{eq:BoxTimesOper}
 (p\otimes a)*_j (q\otimes b)
 =\sum\limits_{i\in I}
  (p\circ^i q) \otimes (a *_{i,j} b),
  \quad p,q\in B,\ a,b\in A.
\end{equation}

\begin{proposition}\label{prop:PropBlackDual}
An algebra $A$ 
with bilinear operations $(\cdot *_{i,j} \cdot)$
is a $(\mathcal P \bullet\mathcal Q)$-algebra
if and only if for every $\mathcal P^!$-algebra $B$
the space $B\otimes A$ equipped with the operations
\eqref{eq:BoxTimesOper} is a $\mathcal Q$-algebra.
\end{proposition}

\begin{proof}
We keep the notations for $\mathcal P$, $\mathcal Q$ introduced above.
Denote by $E$ the subset of $\mathcal F_{V}(3)^\vee\otimes \mathcal F_{W}(3)^\vee $ naturally isomorphic to
\[
\mathop{\fam 0 Ind}\nolimits_{S_2}^{S_3} 
\big (  
(V^\vee \otimes W^\vee) \vecotimes (V^\vee \otimes W^\vee ) 
\big).
\]
Let $\pi $ stand for the conjugate map to the 
embedding $E\hookrightarrow \mathcal F_{V}(3)^\vee\otimes \mathcal F_{W}(3)^\vee $: this is a surjection 
\[
\pi : \big ( \mathcal F_{V}(3)^\vee\otimes \mathcal F_{W}(3)^\vee \big)^\vee\cong \Bbbk_-\otimes \mathcal F_V(3)\otimes \mathcal F_W(3)
\to E^\vee .
\]
By definition, the Manin black product of $\mathcal P$ and $\mathcal Q$ 
is a binary quadratic operad generated by the $S_2$-space
$\Bbbk_-\otimes V\otimes W$
with respect to the relations 
\[
\pi (E^\perp + \Bbbk_-\otimes R\otimes S)\subseteq \Bbbk_-\otimes \mathcal F_V(3)\otimes \mathcal F_W(3).
\]

Denote $e_{ij} = 1_-\otimes e_i\otimes f_j$,
for $i\in I$, $j\in J$.
An arbitrary element $h\in S\subseteq \mathcal F_W(3)$
may be presented as
\[
h = \sum\limits_{j,k\in J} \big(
\alpha _{jk}(f_j\vecotimes f_k)
+ \beta_{jk}(123)(f_j\vecotimes f_k)
+ \gamma_{jk}(132)(f_j\vecotimes f_k)
\big).
\]
Here we choose $e,(123),(132)$ as a basis of $\Bbbk S_3$
over $\Bbbk S_2$.

For every $\mathcal P^!$-algebra $B$ the operation 
$f_j\vecotimes f_k$ applied to 
$p_1\otimes a_1,p_2\otimes a_2,p_3\otimes a_3\in B\otimes A$ takes the same value as
\[
\sum\limits_{i,l\in I} 
(e_i^\vee \vecotimes e_l^\vee)\otimes (e_{ij}\vecotimes e_{lk})
\]
applied to $p_1\otimes p_2\otimes p_3\otimes a_1\otimes a_2\otimes a_3$.
Note that 
\begin{multline}\label{eq:EachSum-Black}
\sum\limits_{i,l\in I} 
(e_i^\vee \vecotimes e_l^\vee)\otimes (e_{ij}\vecotimes e_{lk})
=
\sum\limits_{i,l\in I} 
(e_i^\vee \vecotimes e_l^\vee)\otimes 
((1_-\otimes e_i\otimes f_j)\vecotimes (1_-\otimes e_l\otimes f_k)) \\
=
\sum\limits_{i,l\in I} 
(e_i \vecotimes e_l)^\vee
\otimes 
\pi (1_-\otimes (e_i\vecotimes e_l)\otimes (f_j\vecotimes f_k)) \\
=
(\mathrm {id}_{(V\vecotimes V)^\vee} \otimes \pi)
\sum\limits_{i,l\in I} 
(e_i \vecotimes e_l)^\vee
\otimes 
1_-\otimes (e_i\vecotimes e_l)\otimes (f_j\vecotimes f_k) \\
=
(\mathrm {id}_{\mathcal F_V(3)^\vee} \otimes \pi)
\sum\limits_{i,l\in I}
(e+(123)+(132))
((e_i \vecotimes e_l)^\vee\otimes 1_-\otimes (e_i\vecotimes e_l))\otimes (f_j\vecotimes f_k)
\end{multline}
For the last equation, we added elements with permutations $(123)$ and $(132)$ 
from the kernel of $\pi $.

Therefore, $h(p_1\otimes a_1,p_2\otimes a_2,p_3\otimes a_3)$ may be presented as 
\[
H = \sum\limits_{\xi } F_\xi \otimes G_\xi ,
\quad 
F_\xi \in \mathcal F_{V^\vee}(3), 
\ 
G_\xi \in \mathcal F_{1_-\otimes V\otimes W}(3),
\]
acting on the corresponding element of $B^{\otimes 3}\otimes A^{\otimes 3}$.
It follows from \eqref{eq:EachSum-Black} that 
\[
H = (\mathrm {id}_{\mathcal F_V(3)^\vee} \otimes \pi)\big( 
C^{13}\otimes h
\big ) .
\]
Here $C^{13}$ is the Casimir element of the space $\mathcal F_V(3)$, i.e.,  
\[
C^{13} = \sum\limits_{\mu} \mu^\vee \otimes 1_-\otimes \mu,
\]
the summation is made over a basis of $\mathcal F_V(3)$.

Since the operations $\mathrm{id}\otimes \pi $ and 
$\langle f,\cdot \rangle \otimes \mathrm{id}$ (for $f\in \mathcal F_V(3)$)
commute, we have
\begin{equation}\label{eq:MainIdentity}
(\langle f,\cdot \rangle \otimes \mathrm{id} )H 
= \pi((\langle f,\cdot \rangle \otimes \mathrm{id} )(C^{13}\otimes h))
=\pi (1_-\otimes f\otimes h).
\end{equation}

Suppose $A$ is a $\mathcal P\bullet \mathcal Q$-algebra.
Then for every $f\in R$ the right-hand side of \eqref{eq:MainIdentity}
is a zero map, so $h(p_1\otimes a_1,p_2\otimes a_2, p_3\otimes a_3)=0$
for every $h\in S$. Hence, $B\otimes A$ is a $\mathcal Q$-algebra.

Conversely, assume $B\otimes A$ is a $\mathcal Q$-algebra for every 
$\mathcal P^!$-algebra $B$. In particular, it is true for the free 
$\mathcal P^!$-algebra generated by $p_1,p_2,p_3$. In this case, 
the left-hand side of \eqref{eq:MainIdentity} is a zero map for every
$f\in R$, $h\in S$, so the right-hand side is an identity on $A$, i.e., 
$A$ is a $\mathcal P\bullet \mathcal Q$-algebra.
\end{proof}

\begin{example}
Let $\mathcal P=\di\Lie$ be the operad of Leibniz algebras, and let $\mathcal Q=\Nov$. Let us calculate 
$\di\Lie\bullet \Nov$ by means of Proposition~\ref{prop:PropBlackDual}.
\end{example}

Assume $\di\Lie(2)$ is spanned by $e_1$, $e_2=(12)e_1$, and $\Nov(2)$ is spanned by $f_1$, $f_2=(12)f_1$. Then $(\di\Lie)^!=\pre\Com$, i.e., the operad of Zinbiel algebras satisfying 
the identity
\[
(x_1x_2)x_3 = x_1(x_2x_3) + x_1(x_3x_2)
\] 
for $x_1x_2=e_1^\vee $ (then $e_2^\vee = -x_2x_1$). 

In order to get a result in terms of familiar binary operations, 
$e_{11} = 1_-\otimes e_1\otimes f_1$ by 
$x_1\dashv x_2$ and $e_{21}=1_-\otimes e_2\otimes f_1$
by $-x_1\vdash x_2$,
then 
$e_{12}= -(12)e_{21} = x_2\vdash x_1$,
$e_{22} = -(12)e_{11} = -x_2\dashv x_1$.

Suppose $A$ is an algebra with operations 
$a\vdash b$, $a\dashv b$. Then for every Zinbiel algebra $B$ we need 
$B\otimes A$ with operation 
\[
e_1^\vee \otimes e_{11} + e_2^\vee \otimes e_{21}
\]
to be a Novikov algebra. According to the fixed notations, 
the latter operation is given by
\begin{equation}\label{eq:ZinbNov}
(x\otimes a)(y\otimes b) = xy\otimes a\dashv b + yx\otimes a\vdash b,
\quad x,y\in B,\ a,b\in A.
\end{equation}
The identities on $\vdash$, $\dashv $ should make 
the operation \eqref{eq:ZinbNov} left-symmetric and right commutative,
see Example~\ref{exmp:Nov}.

From the right-commutative identity, we have the equality of
\begin{multline*}
((x\otimes a)(y\otimes b))(z\otimes c) 
\\
=(xy)z\otimes (a\dashv b)\dashv c
+z(xy)\otimes (a\dashv b)\vdash c  
+(yx)z\otimes (a\vdash b)\dashv c
+z(yx)\otimes (a\vdash b)\vdash c 
\\
=(x(yz)+x(zy))\otimes (a\dashv b)\dashv c
+z(xy)\otimes (a\dashv b)\vdash c  \\
+(y(xz)+y(zx))\otimes (a\vdash b)\dashv c
+z(yx)\otimes (a\vdash b)\vdash c 
\end{multline*}
and
\begin{multline*}
((x\otimes a)(z\otimes c))(y\otimes b) 
=(x(zy)+x(yz))\otimes (a\dashv c)\dashv b
+y(xz)\otimes (a\dashv c)\vdash b 
\\
+(z(xy)+z(yx))\otimes (a\vdash c)\dashv b
+y(zx)\otimes (a\vdash c)\vdash b. 
\end{multline*}
Compare similar terms with (linearly independent) monomials 
$x(yz)$, $y(xz)$, $z(yx)$, etc., in the first tensor factor, 
 we obtain
\begin{gather}
(a\dashv b)\dashv c=(a\dashv c)\dashv b,\label{eq:RComRep-1}\\
(a\dashv b)\vdash c = (a\vdash c)\dashv b,
\label{eq:RComRep-Zero}\\
(a\vdash b)\dashv c = (a\vdash c)\vdash b.
 \label{eq:RComRep-2}
\end{gather}
The relation \eqref{eq:RComRep-Zero} may be replaced with 
\begin{equation}\label{eq:RComRep-01}
    (a\dashv b)\vdash c = (a\vdash b)\vdash c
\end{equation}
modulo \eqref{eq:RComRep-2}.

In a similar way, the identity of left symmetry 
\begin{multline*}
((x\otimes a)(y\otimes b))(z\otimes c)-(x\otimes a)((y\otimes b)(z\otimes c))=\\
((y\otimes b)(x\otimes a))(z\otimes c)-(y\otimes b)((x\otimes a)(z\otimes c)).
\end{multline*}
for 
\eqref{eq:ZinbNov} implies the following relations on $\vdash$ and $\dashv$:
\begin{gather}
(a\dashv b)\dashv c - (b\vdash a)\dashv c
= a\dashv (b\dashv c) - b\vdash (a\dashv c), 
\label{eq:LSymRep-1} \\
(a\dashv b)\dashv c - (b\vdash a)\dashv c
= a\dashv (b\vdash c) - b\vdash (a\dashv c), 
\label{eq:LSymRep-Zero1} \\
(a\dashv b)\vdash c - (b\vdash a)\vdash c
= a\vdash (b\vdash c) - b\vdash (a\vdash c), 
\label{eq:LSymRep-Zero2} \\
(a\vdash b)\vdash c - (b\dashv a)\vdash c
= a\vdash (b\vdash c) - b\vdash (a\vdash c). 
\label{eq:LSymRep-2} 
\end{gather}
It follows from \eqref{eq:LSymRep-Zero1} and \eqref{eq:LSymRep-1} that 
\begin{equation}\label{eq:RComRep-02}
    a\dashv (b\vdash c ) = a\dashv (b\dashv c).
\end{equation}
Therefore, the complete list of independent identities that 
hold on $A$ coincides with the replication of right commutativity and left symmetry, i.e., 
\[
\di\Lie \bullet \Nov =\di\Nov.
\]

\begin{theorem}\label{thm:BlackDong}
Given two Dong operads $\mathcal P$ 
and $\mathcal Q$, their Manin black product 
$\mathcal P\bullet \mathcal Q$ also 
satisfies the Dong Property. 
\end{theorem}

\begin{proof}
Let us follow the notations above: the operations
on a $\mathcal P^!$-algebra
are denoted $(\cdot \circ^i \cdot)$, $i\in I$,
the operations on a $\mathcal Q$-algebra are
$(\cdot *_j \cdot)$, $j\in J$, and let 
$(\cdot *_{i,j} \cdot )$, $i\in I$, $j\in J$, 
stand for the operations on a 
$(\mathcal P\bullet\mathcal Q)$-algebra.

Suppose we are given three pairwise local 
formal distributions 
$a(z),b(z),c(z)$ over a 
$(\mathcal P\bullet\mathcal Q)$-algebra $A$.
Recall that 
\[
x(z) = \sum\limits_{s\in \mathbb Z} x(s)z^{-s-1}, \quad 
x\in \{a,b,c\}, \ x(s)\in A.
\]
Consider the free $\mathcal P^!$-algebra $B$ generated 
by three elements $p,q,t$.
Then 
\[
\bar a(z)=p\otimes a(z),\quad 
\bar b(z)=q\otimes b(z),\quad 
\bar c(z)=t\otimes c(z)
\]
are formal distributions over the 
$\mathcal Q$-algebra $B\otimes A$. 
It is obvious that $\bar a(z),\bar b(z), \bar c(z)$
are pairwise local (with respect to the operations 
$*_j$, $j\in J$).
Since $\mathcal Q$ is a Dong operad, for every $n\in \mathbb Z_+$ and for every $j\in J$
the distributions 
\[
(\bar a\ooc{j}{n} \bar b)(z),\, \bar c(z) 
\in (B\otimes A)[[z,z^{-1}]]
\]
are local. The latter means that
for every $j_1\in J$ there exists $N\ge 0$
such that 
\[
(w-z)^N(\bar a\ooc{j}{n} \bar b)(w) *_{j_1} \bar c(z) =0.
\]
By definition, see \eqref{eq:n-prod},
we have 
\[
(\bar a\ooc{j}{n} \bar b)(w)
=\Res\limits_{\xi=0} (\xi-w)^n \bar a(\xi)*_j \bar b(w).
\]
It follows from \eqref{eq:BoxTimesOper} that
\[
\bar a(\xi)*_j \bar b(w)
=\sum\limits_{i\in I} (p\circ ^i q)\otimes (a(\xi) *_{j,i} b(w) ), 
\]
hence, 
\begin{multline*}
0 = (w-z)^N(\bar a\ooc{j}{n} \bar b)(w) *_j \bar c(z) 
\\
= (w-z)^N \Res\limits_{\xi=0} (\xi-w)^n
\bigg (\sum\limits_{i\in I} (p\circ ^i q)\otimes (a(\xi) *_{i,j} b(w) ) \bigg) *_{j_1} (t\otimes c(z))
\\
=
 (w-z)^N 
\sum\limits_{i,i_1\in I} (p\circ ^i q)\circ^{i_1} t
\otimes \Res\limits_{\xi=0}(\xi-w)^n(a(\xi) *_{i,j} b(w) ) *_{j_1,i_1} c(z)
\\
=
\sum\limits_{i,i_1\in I}
(p\circ^i q)\circ^{i_1} t
\otimes 
(w-z)^N (a\ooc{j,i}{n} b)(w) *_{i_1,j_1} c(z).
\end{multline*}
Since the elements $(p\circ^i q)\circ^{i_1} t$
are linearly independent in $B$ by Theorem~\ref{thm:Necessary}, 
we have 
\[
(w-z)^N (a\ooc{i,j}{n} b)(w) *_{i_1,j_1} c(z) =0\in 
A[[z,z^{-1},w,w^{-1}]]
\]
for every $i,i_1\in I$,
as desired.
\end{proof}

\subsection*{Acknowledgments}
We are grateful to Vladimir Dotsenko 
for the idea to study the analogues of the Dong Lemma
at the level of operads. This research was funded by the Science Committee of the Ministry of Science and Higher Education of the Republic of Kazakhstan (Grant No. AP26193761).
The first author was supported by the Program of Fundamental Research RAS (project FWNF-2022-0002).

%
%


\begin{thebibliography}{99}

\bibitem{BaiKo-TPalgebras}
Bai C., Bai R., Guo L., Wu Y., Transposed Poisson algebras, Novikov-Poisson algebras and 3-Lie algebras, Journal of Algebra, 2023, 632, 535--566.

\bibitem{BBGN13}
Bai C., Bellier O., Guo L., Ni X., Splitting of operations, Manin products, and Rota-Baxter operators, International Mathematics Research Notices, 2013, 2013(3), 485–524.

\bibitem{BDK2001}
Bakalov B., D'Andrea A., Kac V.G., Theory of finite pseudoalgebras, Advances in Mathematics, 2001, 162(1), 1–140.

\bibitem{BK-et-al2018}
Bakalov B., De Sole A., Heluani R., Kac V. G., An operadic approach to vertex algebra and Poisson vertex algebra cohomology, Japanese Journal of Mathematics, 2019, 14(2), 249–342.

\bibitem{BalNov}
Balinskii A. A., Novikov S. P., Poisson brackets of hydrodynamic type, Frobenius algebras and Lie algebras, Sov. Math. Dokl., 32(1) (1985), 228--231.

\bibitem{BPZ198x}
Belavin A.A., Polyakov A.M., Zamolodchikov A.B., Infinite conformal symmetry in two-dimensional quantum field theory, Nuclear Physics, Section B, 1984, 241(2), 333–380.

\bibitem{BokKol2024}
Bokut L. A., Kolesnikov P. S. On the locality of formal distributions over right-symmetric and Novikov algebras, The Bulletin of Irkutsk State University. Series Mathematics, 2024, vol. 50, pp. 83–100. (in Russian). ArXiv:2404.17232v1.

\bibitem{Bor198x}
Borcherds R. E., Vertex algebras, Kac-Moody algebras, and the monster, Proceedings of the National Academic Society USA, 83 (1986), 3068–3071.

\bibitem{BrD_Companion}
Bremner M. R., Dotsenko  V., Algebraic Operads An Algorithmic Companion, Chapman Hall, 2016.

\bibitem{DzhLofwall2002}
Dzhumadil'daev A., L\"ofwall C., Trees, free right-symmetric algebras, free Novikov algebras and identities, Homology, Homotopy and Applications, 
2002, 4(2 I), 165--190.

\bibitem{FbZ2003}
Frenkel E., Ben-Zvi D., Vertex algebras and algebraic curves, American Mathematical Society, 88, 2004.

\bibitem{GelDorfm}
Gelfand I. M., Dorfman I. Ya., Hamilton operators and associated algebraic structures, Functional analysis and its application, 13 (1979), no. 4, 13--30.

\bibitem{GinKapr1994}
Ginzburg V., Kapranov M., Koszul duality for operads, Duke Mathematical Journal, 76(1) (1994), 203--272.

\bibitem{GK13}
Gubarev V. Y., Kolesnikov P. S., Embedding of dendriform algebras into Rota-Baxter algebras, Central European Journal of Mathematics, 2013, 11(2), 226–245.

\bibitem{GK14}
Gubarev V. Yu., Kolesnikov P. S.,
Operads of decorated trees and their duals,
Comment. Math. Univ. Carolin. 55(4) (2014) 421--445.

\bibitem{Hong2015}
Hong Y., Li F., Left-symmetric conformal algebras and vertex algebras, Journal of Pure and Applied Algebra, 2015, 219(8), 3543–3567.

\bibitem{KacVABeginn}
Kac V. G., Vertex algebras for beginners, American Mathematical Society 10, 1998.

\bibitem{KacRetakh-JordSuper}
Kac V. G., Retakh A.,
Simple Jordan conformal superalgebras,
J. Algebra Appl. 
7(4) (2008) 517--533.

\bibitem{Kol2006CA}
Kolesnikov P., Identities of conformal algebras and pseudoalgebras, Communications in Algebra, 2006, 34(6), 1965--1979.

\bibitem{HLie199x}
Li H., 
Local systems of vertex operators, vertex superalgebras and modules, 
J. Pure Appl. Alg. 109 (1996) 143--195.

\bibitem{LodayPirash199?}
Loday J.-L., Pirashvili T.,
Universal envelopping algebras of Leibniz algebras
and homology,
Math. Ann. 296 (1993) 139--158.

\bibitem{LodayVallette_book}
Loday J. L., Vallette B., Algebraic operads, Springer Science and Business Media, 2012.

\bibitem{MarklRemm}
Markl M., Remm E., Algebras with one operation including Poisson and other Lie-admissible algebras, Journal of Algebra, 2006, 299(1), 171--189.

\bibitem{Osborn90}
Osborn J. M., Novikov algebras, Nova J. Algebra \& Geom., 1 (1992), 1--14.

\bibitem{BGG_HK_Proceed}
Pei J., Bai C., Guo L., Ni X., 
Replicators, Manin white product of binary operads 
and average operators. 
In: 
New trends in algebras and combinatorics, pp. 317--353. World Scientific Publishing Co. Pte. Ltd.,
Hackensack, NJ (2020).

\bibitem{Roit1999}
Roitman M., On free conformal and vertex algebras, Journal of Algebra, 1999, 217(2), 496--527.

\bibitem{Wightman5x}
Wightman A. S.,
Quantum field theory in terms of vacuum expectation values. 
Phys. Rev. 101 (1956), 860--866.

\bibitem{Xu1993}
Xu X., On simple Novikov algebras and their irreducible modules, Journal of Algebra, 1996, 185(3), 905--934.

\bibitem{Xu1999}
Xu X., Quadratic conformal superalgebras, Journal of Algebra, 2000, 231(1), 1--38.

\bibitem{Yuan22}
Yuan L., $O$-operators and Nijenhuis operators of associative conformal algebras,
J. Algebra 609 (2022) 245--291.

\bibitem{Zelmanov87}
Zelmanov E. I., On a class of local translation invariant Lie algebras, English transi. Soviet Math. Dokl. 35 (1987), 216--218.

\end{thebibliography}
\end{document}